\documentclass[11pt,a4paper,leqno,twoside, headinclude]{amsart}

\usepackage[tracking=true]{microtype}
\DeclareMicrotypeSet*[tracking]{my}%
  { font = */*/*/sc/* }%
\SetTracking{ encoding = *, shape = sc }{ 45 }

\title[Homology vs. homotopy]{Homology versus homotopy in fibrations\\ and in limits}
\def\titl{}
\def\auth{Manuel Amann}
\date{June 1st, 2020}

\thanks{}

\author{\auth}

\usepackage[english]{babel}

\usepackage[colorlinks,pdftex, plainpages=false]{hyperref}
\hypersetup{pdftitle=\titl, pdfauthor=\auth, pdftoolbar=false,
plainpages=false, hyperindex=true, pdfdisplaydoctitle=true}

\hypersetup{colorlinks,%
citecolor=black,%
filecolor=black,%
linkcolor=black,%
urlcolor=black,%
pdftex}

\usepackage{color}

\usepackage{amsmath, amssymb, amscd, amsthm}
\usepackage{stmaryrd}
\usepackage{latexsym}
\usepackage[all]{xy}
\usepackage{pb-diagram}
\usepackage{rotating}
\usepackage{multicol}
\usepackage{lscape}
\usepackage{wasysym}
\usepackage{longtable}
\usepackage{enumerate}

\xyoption{all}

\newtheorem{theo}{Theorem}[section]
\newtheorem{main}{Theorem}
\newtheorem{maincor}[main]{Corollary}
\newtheorem*{main*}{Theorem}

\newtheorem*{mainprop*}{Proposition}

\newtheorem{mainconj}{Conjecture}

\newtheorem{prop}[theo]{Proposition}
\newtheorem{defi2}[theo]{Definition}
\newtheorem*{defi2*}{Definition}
\newenvironment{defi}{\begin{defi2}\normalfont}{\end{defi2}}
\newenvironment{defi*}{\begin{defi2*}\normalfont}{\end{defi2*}}

\newenvironment{defin*}[1]{\begin{defi2*}[#1]\normalfont}{\end{defi2*}}
\newtheorem*{rem2*}{Remark}
\newenvironment{rem*}{\begin{rem2*}\normalfont}{\hfill$\boxbox$\end{rem2*}}
\newtheorem{rem2}[theo]{Remark}
\newenvironment{rem}{\begin{rem2}\normalfont}{\hfill$\boxbox$\end{rem2}}

\newtheorem{lemma}[theo]{Lemma}
\newtheorem{cor}[theo]{Corollary}
\newtheorem*{cor*}{Corollary}

\newtheorem{conj}[theo]{Conjecture}

\newtheorem*{conj*}{Conjecture}
\newtheorem*{theo*}{Theorem}
\newtheorem*{ques*}{Question}
\newtheorem*{mi2}{Main Idea}

\newtheorem{ex2}[theo]{Example}
\newenvironment{ex}{\begin{ex2}\normalfont}{\hfill$\boxbox$\end{ex2}}

\newtheorem{exer2}[theo]{Exercise}

\newtheorem{alg2}[theo]{Algorithm}

\newcommand{\cc}{{\mathbb{C}}}                                     
\newcommand{\hh}{{\mathbb{H}}}                                     
\newcommand{\nn}{{\mathbb{N}}}                                     
\newcommand{\qq}{{\mathbb{Q}}}                                     
\newcommand{\rr}{{\mathbb{R}}}                                     
\newcommand{\pp}{{\mathbf{P}}}                                     
\newcommand{\s}{{\mathbb{S}}}                                      
\newcommand{\SU}{{\mathbf{SU}}}                                    
\newcommand{\Sp}{{\mathbf{Sp}}}                                    
\newcommand{\dif} {{\operatorname{d}}}                             
\newcommand{\In} {{\,\subseteq\,}}                                 
\newcommand{\im} {{\operatorname{im\,}}}                           
\newcommand{\rk}{{\operatorname{rk\,}}}                            
\newcommand{\co}{\colon\thinspace}                                 
\newcommand{\vproof}{{\begin{flushright} \qed                      
                      \end{flushright}}}

\newcommand{\comment}[1]{}                                         
\newcommand{\xto}[1]{\xrightarrow{#1}}                             
\newcommand{\hto}[1]{\overset{#1}{\hookrightarrow}}                

\newcommand{\biq}[2]{#1\;\!\!\!\sslash \;\!\!\!#2}                 
\newcommand{\case}[1]{\textbf{Case #1.}}                           
\newcommand{\ack}{\noindent\textbf{Acknowledgements. }}            
\newcommand{\str}{\noindent\textbf{Structure of the article. }}    
\newcommand{\odd}{\textrm{odd}}                                    
\newcommand{\even}{\textrm{even}}                                  

\newenvironment{prf}{\begin{proof}[\textsc{Proof}]} {\end{proof}}     

\begin{document}

\maketitle \thispagestyle{empty}


\begin{abstract}
Motivated by prominent problems like the Hilali conjecture Yamaguchi--Yokura recently proposed certain estimates on the relations of the dimensions of rational homotopy and rational cohomology groups of fibre, base and total spaces in a fibration of rationally elliptic spaces.

In this article we prove these estimates in the category of formal elliptic spaces and, in general, whenever the total space in addition has positive Euler characteristic or has the rational homotopy type of a homogeneous manifold (respectively of a known example) of positive sectional curvature. Additionally, we provide general estimates approximating the conjectured ones.

Moreover, we suggest to study families of rationally elliptic spaces under certain asymptotics, and we discuss the conjectured estimates from this perspective for two-stage spaces.
\end{abstract}


\section*{Introduction}

The Hilali conjecture (\cite{Hil80}) speculates that for a rationally elliptic space, i.e.~a simply-connected space with both finite rational cohomology and rational homotopy groups, the dimension of rational cohomology is at least as large as the dimension of the rational homotopy groups, i.e.~in other words, their well-defined quotient
\begin{align*}
h(X)= \frac{\dim \pi_*(X)\otimes \qq}{\dim H^*(X)}
\end{align*}
is well-defined and smaller equal one. While the conjecture still being open, this quotient was considered in several different further circumstances (for example see \cite{NY11}).

Recently, it was asked by Yamaguchi--Yokura how this quotient behaves in fibrations of rationally elliptic spaces.

It is the goal of this article to provide several special cases of their conjectured estimates on the one hand and, on the other hand, to study this quotient asymptotically---first suggesting, specifying and discussing different reasonable notions of ``asymptotic behaviour'' for families of rationally elliptic spaces.

\bigskip

Throughout this article we denote by $X$ a simply-connected CW-complex. Cohomology is considered with rational coefficients.
As stated above, we call $X$ \emph{rationally elliptic} if it is simply-connected and both $\dim \pi_*(X)\otimes \qq<\infty$ and $\dim H^*(X)<\infty$. It is called \emph{$F_0$} or \emph{positively elliptic}, if $\chi(X)>0$. (By abuse of notation, we shall refer to rationally elliptic spaces as just being \emph{elliptic}. This is also indebted to the fact that the article builds on rational methods. In particular, whenever we speak of a ``fibration'', it is actually enough to have a ``rational fibration'' structure.)

Recall that the prominent subclass of elliptic spaces, the class of \emph{two-stage spaces}, is defined as follows: Their minimal Sullivan models $(\Lambda V,\dif)$ up to isomorphism admit decompositions of the form $V^\odd=W_0\oplus W_1$ and $\dif(V^\even\oplus W_0)=0$, $\dif W^1\In \Lambda (V^\even \oplus W^0)$.

We provide several notions of convergence for families of elliptic spaces. See Section \ref{sec01} for an elaborate discussion of this. In particular, we rigorously define \emph{$\pi$-convergence} there. It appears to be very hard to control the possible values of $h(X)$; so it seems reasonable to consider its asymptotic behaviour. To our knowledge this is the first time such a discussion is launched. 
As a first step this then permits to prove
\begin{main}\label{theoA}
The family of two-stage spaces $\mathcal{X}$  $\pi$-converges to $0$, i.e.~with $\dim \pi_*(X)\otimes \qq$ (for $X\in \mathcal{X}$) tending to $\infty$, $h(X)$ tends to $0$.
\end{main}

The following results deal with the behaviour of $h(X)$ in fibrations. Hence consider a fibration $F\hto{} X \to B$ of rationally elliptic spaces.
\begin{conj}[Yamaguchi--Yokura, \cite{YY20}]\label{conj01}
\begin{align}\label{eqn02}
\frac{1}{2}\cdot h(F\times B) \leq h(X) < h(F)+h(B)+\frac{1}{4}
\end{align}
\end{conj}

As an application of Theorem \ref{theoA} we discuss an asymptotic version of this problem first.
We say that a class $\mathcal{X}$ of rationally elliptic spaces \emph{asymptotically satisfies} Conjecture \ref{conj01} if the following holds: There is $k\in \nn$ such that if $\dim \pi_*(X)\otimes \qq\geq k$ for $X\in \mathcal{X}$, then $X$ satisfies the conjecture.

Clearly, any family $\mathcal{X}$ $\pi$-converging to $0$ asymptotically satisfies the right hand side equation, i.e.~for $X\in\mathcal{X}$ with large enough rational homotopy groups, the quotient $h(X)$ is smaller than $1/4$.
\begin{maincor}\label{corA}
Let $\mathcal{X}$ be a family $\pi$-convergent to $0$. Then $\mathcal{X}$ asymptotically satisfies the right hand side of Inequality \eqref{eqn02}. In particular, this holds true for two-stage spaces.
\end{maincor}
Clearly, there are two-stage spaces (already products of spheres) of arbitrarily large rational homotopy. Actually, our estimates for two-stage spaces are explicit, and it is easy from there to provide concrete numbers for $\dim \pi_*(X)\otimes \qq$ starting the range in which the inequality holds.

\bigskip

The following results aim to verify Conjecture \ref{conj01} in particular cases. As a vital tool to do so we first verify a more or less close approximation to the left hand side of Conjecture \ref{conj01} in
\begin{main}\label{theoB}
For any fibration $F\hto{}X\to B$ of elliptic spaces it holds that
\begin{align*}
h(F\times B) \leq 3\cdot h(X)
\end{align*}
\end{main}
Next we prove the conjecture whenever $X$ is positively elliptic, as well as in the category of formal elliptic spaces. Recall that a space is \emph{formal} if its rational homotopy type is isomorphic to the one of its cohomology algebra. Formal spaces form one of the most prominent classes of spaces in rational homotopy theory.
\begin{main}\label{theoC}
Let $F\hto{}X\to B$ be a fibration of elliptic spaces.

If $F$ is an $F_0$-space it holds that
\begin{align*}
h(F\times B) \leq 2\cdot h(X)
\end{align*}

Moreover, Conjecture \ref{conj01} holds with respect to any such fibration whenever
\begin{itemize}
\item
$X$ is an $F_0$-space, or
\item
$F$ is an $F_0$-space and satisfies the Halperin conjecture.
\end{itemize}
\end{main}
For a brief discussion of the Halperin conjecture see Section \ref{secrat}. In particular, there we provide a list of several classes of spaces for which the conjecture is verified.

The last formulation is stricter than necessary: For a fixed totally non-homologous to zero fibration $F\hto{} X\to B$ with $F$ an $F_0$-space the required inequalities hold already. With the presented formulation it is our goal to stress that conjecturally any $F_0$-space $F$ should render the fibration totally non-homologous to zero.

Combining, refining and extending the previous arguments we finally obtain
\begin{main}\label{theoD}
Conjecture \ref{conj01} holds for a fibration $F\hto{} X\to B$ of elliptic formal spaces.
\end{main}
This is proved in Propositions \ref{prop01} and \ref{prop02}.

As a corollary to this we can prove Conjecture \ref{conj01} whenever $X$ is a known example of positive sectional curvature, respectively a homogeneous space of positive curvature---see Section \ref{secpos} for more details on these classes. This is particularly interesting for different reasons: First of all these spaces constitute a nice class of highly important geometric examples. Second, maybe more strikingly, let us recall the \emph{Petersen--Wilhelm conjecture} which states that whenever $X\to B$ is a Riemannian submersion with $X$ (and consequently $B$) positively curved Riemannian manifolds, then $2\dim B> \dim X$; respectively, in the case of compact spaces, when this submersion is a fibration $F\hto{}X\to B$, $\dim B> \dim F$. We recall the general property (for example see \cite{AK15}, \cite[Proposition 1, p.~5]{GR18}) that with $X$ being elliptic (and not necessarily a manifold with curvature bound) the fibration features in the category of elliptic spaces already. In particular, this would provide an a priori weaker formulation of Conjecture \ref{conj01}.

In \cite{AK15} we proved the Petersen--Wilhelm conjecture in a much more general context for the known examples of positive curvature of even dimensions only using their rational structure. (Since several odd-dimensional examples rationally split as a product, rational tools are not enough in odd dimensions; in \cite{GR18} the odd-dimensional examples are verified using finite coefficients.)

Note further that due to the Bott--Grove--Halperin conjecture positively curved manifolds should be elliptic. In even dimensions the equally famous Hopf conjecture speculates that they have positive Euler characteristic. That is, conjecturally the case of even-dimensional positively curved manifolds should be completely covered by Theorem \ref{theoB}.

In summary, in this context Conjecture \ref{conj01} controls much more complicated invariants of fibration decompositions of positively curved manifolds than merely dimensions. Given all previous observations and conjectures one might speculate that

\begin{conj*}
For a closed simply-connected positively curved manifold $M$ and any fibration $F\hto{}M\to B$ of simply-connected spaces Estimate \eqref{eqn02} is well-defined and holds true.
\end{conj*}

Viewed from a different angle, we once again observe that positively curved manifolds seem to constitute a class of spaces behaving extremely well with respect to several different topological approaches. We verify this speculation on the known examples.
\begin{maincor}\label{corB}
Conjecture \ref{conj01} holds true whenever the cohomology algebra $H^*(X)$ is generated by at most one even-degree and at most one odd-degree element. In particular, this is true if $X$ has the rational type of a simply-connected closed homogeneous space of positive sectional curvature respectively of any known example of a closed manifold admitting positive sectional curvature.
\end{maincor}
The content of this is the observation that any fibration then only involves formal spaces.
We remark that the confirmation of the conjecture for $X$ positively elliptic is yet another corollary of Theorem \ref{theoD} as well: If $\chi(X)>0$, by the multiplicativity of the Euler characteristic in fibrations, so are $\chi(B),\chi(F)>0$. It is well-known that $F_0$-spaces are formal.

\bigskip

We leave it to the reader to reformulate our results in larger generality for nilpotent spaces and nilpotent fibrations.

\bigskip

\str In Section \ref{sec09} we discuss some relevant aspects from Rational Homotopy Theory. In Section \ref{sec01} we note first observations on the conjecture before, in the second part of the section, we elaborately consider and discuss different notions of convergence for families of elliptic spaces. In Section \ref{secA} we explain the proof of Theorem A, which is rather independent of the following arguments. Section \ref{secB} is devoted to the proof of Theorem \ref{theoB}. In particular, there we prove Lemma \ref{lemma01}, which is central to our arguments and underlies nearly all further (and partly even previous) work. Finally, in Section \ref{secCD} we refine and massively extend the previous arguments in order to provide proofs of Theorems \ref{theoC} and \ref{theoD}. As an application of the obtained results we use this to show Corollary \ref{corB} in a subsequent step.

\bigskip

\ack The author was supported both by a Heisenberg grant and his research grant AM 342/4-1 of the German Research Foundation; he is moreover associated to the DFG Priority Programme 2026.


\section{Some tools from Rational Homotopy Theory}\label{sec09}

\subsection{Excerpts from Rational Homotopy Theory}\label{secrat}

This section cannot provide and is not intended to give an introduction to the theory. We expect the reader to have gained a certain familiarity with necessary concepts for example from \cite{FHT01} or \cite{FOT08}. We merely recall some tools and aspects which play a larger role in the article.

Many computations of the article rely on the theory of \emph{(minimal) Sullivan models} of simply-connected spaces $X$. Just to recall these are certain commutative differential graded algebras $(\Lambda V,\dif)$ encoding the rational homotopy type of $X$ with $V$ a positively graded rational vector space and $\Lambda V$ the tensor product on the symmetric algebra of the evenly-graded part $V^\even$ and the exterior algebra on the oddly-graded part $V^\odd$. Moreover, we use \emph{relative models} and \emph{models of fibrations} as constructed in \cite[Proposition 15.5]{FHT01}. That is, for a fibration of simply-connected spaces $F\hto{} E \to B$ and for Sullivan models $(\Lambda V,\bar\dif)$ of $F$ and $(\Lambda W,\dif)$ of $B$, a model for $E$ is given by a tensor product $(\Lambda V\otimes \Lambda W,\dif)$ where $(\Lambda W,\dif)$ is a differential subalgebra, and the projection induced by $W\to 0$ yields the model $(\Lambda V,\bar \dif)$ of $F$.

We investigate fibrations and their Sullivan models from a cohomological and a homotopical point of view. For the associated Serre spectral sequence see \cite[Chapter 18]{FHT01}. For the associated long exact sequence of homotopy groups in terms of models see \cite[Section 15(e), p.~214]{FHT01}. In particular, recall that with the terminology from the last paragraph, the long exact homotopy sequence dualises to the exact sequence
\begin{align*}
\ldots \to W^k\to H^k(W\oplus V,\dif_0) \to V^k \xto{\dif_0} W^{k+1}\to \ldots
\end{align*}
with transgression $\dif_0$ with respect to which also cohomology $H^k(W\oplus V,\dif_0)$ is taken. This $\dif_0$ denotes the linear part of the differential $\dif$ on $(\Lambda (V\oplus W)$ defined by $\im (\dif -\dif_0) \In \Lambda^{\geq 2} (V\oplus W)$. (Clearly, $V^k\cong \pi_k(F)\otimes \qq$, $W^k\cong \pi_k(B)\otimes \qq$---see \cite[Theorem 15.11, p.~208]{FHT01}---and $H^k(W\oplus V,\dif_0)\cong \pi_k(E)\otimes \qq$ taking into account that the model of the fibration is not necessarily minimal.)

We shall speak of rational homotopy groups of $F$ and $B$ being \emph{contracted} when passing to $X$ which is supposed to indicate that such a homotopy group lies in the kernel respectively the image of $\dif_0$ and hence exists in $F$ respectively $B$, but no longer contributes non-trivial homotopy to $X$.

We shall moreover draw on Euler and homotopy Euler characteristics. We use the convention to define the latter for an elliptic minimal Sullivan algebra $(\Lambda V,\dif)$ as
\begin{align*}
\chi_\pi(\Lambda V,\dif)=\dim V^\odd-\dim V^\even
\end{align*}
The Euler characteristic is multiplicative (which can be proved using the Serre spectral sequence), the homotopy Euler characteristic is additive in fibrations (as follows from the depicted long exact homotopy sequence).

\bigskip

The \emph{formal dimension} $n$ of an elliptic space, i.e.~the largest degree with non-trivial cohomology can be computed via the following \emph{dimension formula} using the degrees and dimensions of its homotopy groups---see \cite[Theorems 32.2 (iii), p.~436 and 32.6 (i), p.~441]{FHT01}. For this we recall the \emph{even and odd exponents} $a_i$ and $b_i$ of a minimal Sullivan algebra $(\Lambda V,\dif)$ defined by the property that the $2 a_i$ are the degrees of a basis of $V^\even$ and the $2 b_i-1$ are the degrees of a basis of $V^\odd$. It then holds that
\begin{align*}
\sum (2 b_i-1) -\sum (2a_i-1)=n
\end{align*}
Compare Remark \ref{remdim}.

\bigskip

Elliptic spaces $X$ of positive Euler characteristic, so-called \emph{$F_0$-spaces} or \emph{positively elliptic spaces} possess a very rigid structure: Their rational cohomology is concentrated in even degrees; actually it is given by a polynomial algebra modulo a regular sequence whence these spaces are (hyper-/intrinsically) formal. Moreover, from \cite[Formula (32.14), p.~446]{FHT01} we recall that the total dimension of their cohomology, i.e.~their Euler characteristic, which equals the sum of all Betti numbers in this case, is given by
\begin{align}\label{eqn098}
\dim H^*(X)=\prod_{i=1}^q \frac{b_i}{a_i}
\end{align}
where the $b_i$ and $a_i$ range over the odd respectively the even exponents of a minimal model of $X$.

As a consequence, positively elliptic spaces admit \emph{pure models}---see \cite[Chapter 32, p.~434]{FHT01}. This contributes to the importance of pure spaces in rational homotopy theory. There are many prominent classes of pure spaces featuring homogeneous spaces and biquotients as well as cohomogeneity one manifolds. Recall our definition of \emph{two-stage spaces} in the introduction which clearly constitutes a slight generalisation of pureness. Two-stage spaces gain special importance due to the following
\begin{prop}[see Proposition 5.10, p.~32, in \cite{AZ19}, cf.~\cite{KY05}]\label{prop05}
Let $X$ be a formal elliptic space. Then rationally it is the total space of a totally non-homologous to zero fibration with model
\begin{align*}
(\Lambda B,0)\rightarrow (\Lambda B\otimes \Lambda V,\dif)\rightarrow (\Lambda V,\bar \dif)
\end{align*}
where $B=B^{odd}$, and $(\Lambda V,\bar \dif)$ is positively elliptic.
\end{prop}
That is, in particular, formal elliptic spaces are two-stage.

Recall that a fibration $F\hto{}X \to B$ is called \emph{totally non-homologous to zero} if the induced map $H^*(X)\to H^*(F)$ is surjective, or, equivalently, if the associated Serre spectral sequence degenerates at the $E_2$-term.

\begin{rem}\label{remmin}
Moreover, we may assume that the model $(\Lambda B\otimes \Lambda V,\dif)$ is \emph{minimal}, i.e.~the fibration to be $\pi$-trivial as well, (and then decompose it as depicted). This follows from the proof of \cite[Proposition 5.10]{AZ19} in which we decomposed a two-stage model of $X$ with stage one mapping to a regular sequence in the algebra generated by stage $0$ in this form. Without restriction, we may choose the model we start with to be minimal. Indeed, the model comes from \cite[Theorem II, p.~577]{FH82}, and can be chosen as a minimal model of the hyperformal cohomology $H^*(X)$.
\end{rem}

\bigskip

One of the most famous and most influential conjectures in the area is
\begin{conj}[Halperin]
Let $F\hto{} X\to B$ be a fibration of simply-connected spaces with $F$ positively elliptic. Then the fibration is totally non-homologous to zero.
\end{conj}

In particular, this conjecture was verified
\begin{itemize}
\item
on compact homogeneous spaces of positive Euler characteristic (\cite{ST87}).
\item
on simply-connected Hard-Lefschetz spaces (\cite{Mei83}).
\item in the case of at most three generators of the cohomology algebra (see \cite{Tho81} and \cite{Lup90}).
\item for spaces of formal dimension at most $16$ or Euler characteristic at most $16$ (see \cite[Theorem 11.6]{AK15b}).
\item
in the ``generic case'' (cf.~\cite{PP96}).
\item
to be closed under fibrations of simply-connected spaces of finite type (cf.~\cite{Mar90}), i.e.~if both base and fibre satisfy the Halperin conjecture, so does the total space.
\end{itemize}

These classes of examples enrich Theorem \ref{theoC}.

Note further that it is known that the Halperin conjecture for an elliptic space $F$ holds if and only if it holds in the category of elliptic spaces, more precisely, already for base spaces being odd dimensional spheres (see \cite[Theorem 1.5, p.~6]{Lup98}, \cite{Mei82}).

\subsection{Positively curved spaces and their rational structure} \label{secpos}
By ``positive curvature'' we shall always denote positive \emph{sectional} curvature.

The known examples of simply-connected positively curved closed manifolds are the following (cf.~\cite{GR18}):
\begin{itemize}
\item the subsequent homogeneous spaces, namely compact rank one symmetric spaces $\s^n$, $\cc\pp^n$, $\hh\pp^n$, $\operatorname{CaP}^2$, the Wallach flag manifolds $W^6$, $W^{12}$, $W^{24}$, the Aloff--Wallach spaces $W_{p,q}^7$, and the Berger spaces $B^7$, $B^{13}$.
\item the biquotients $E^6$ due to Eschenburg, the family $\biq{\SU(3)}{\s^1}$ \linebreak[4](parametrised by different inclusions) generalising and comprising the Aloff--Wallach spaces, the family of Bazaikin spaces $\biq{\SU(5)}{\Sp(2)\s^1}$ in dimension $13$ containing  $B^{13}$, and
\item a cohomogeneity one example $P^2$ of dimension $7$ due to Dearricot and Grove--Verdiani--Ziller.
\end{itemize}
Without going into details, collecting the information for example from \cite{Zil12}, \cite{AK15}, \cite{{GR18}} we derive that all these spaces are formal and elliptic, and, in any case, the following holds:
\begin{itemize}
\item If $M$ is even-dimensional, it is positively elliptic.
\item If $M$ is odd-dimensional, it satisfies $\chi_\pi(M)=1$. It is either rationally a sphere, or its rational cohomology algebra has exactly one generator in positive even-degree and exactly one in odd-degree.
\end{itemize}
This is the necessary information underlying the proof of Corollary \ref{corB} for positively curved manifolds (see Property $(*)$ and Remark \ref{rem01}).

We remark further that there is a classification of simply-connected positively curved homogeneous spaces by Wallach and B\'erard--Bergery (for example see \cite{WZ15}) which states that the cited homogeneous examples are actually the only ones.


\section{First observations}\label{sec01}

\subsection{Fibrations}\label{secfib}
Let $F\hto{} X\to B$ a fibration of elliptic spaces.
We call it \emph{$\pi$-trivial}, if $\pi_*(X)\otimes \qq=\pi_*(F)\otimes \qq\oplus \pi_*(B)\otimes \qq$, or, equivalently, if the relative minimal model of the fibration is actually a minimal model of $X$.

It is interesting to observe that $\pi$-trivial fibrations play a role converse to the one of totally non-homologous to zero ones with respect to Conjecture \ref{conj01}; more precisely,
\begin{itemize}
\item if the fibration is totally non-homologous to zero, then $h(X)\leq h(F) + h(B)$ (note that the computations of \cite[Page 3]{YY20} for the product fibration apply similarly to yield this inequality), and the right hand side of \eqref{eqn02} is satisfied, in particular.
\item if the fibration is $\pi$-trivial, then $h(F\times B)\leq h(X)$, and the left hand side of \eqref{eqn02} is satisfied, in particular.
\end{itemize}
As for the latter, it suffices to recall that the Serre spectral sequence of the fibration $F\hto{} X\to B$ of simply-connected spaces of finite-dimensional rational cohomology satisfies $E_2^{p,q}=H^p(B)\otimes H^q(F)$, whence
\begin{align}\label{eqn123}
\dim H^*(X)\leq \dim H^*(F)\cdot \dim H^*(B)
\end{align}
If the fibration is $\pi$-trivial, it follows that $\pi_*(F\times E)\otimes \qq=\pi_*(X)\otimes \qq$. We deduce the given estimate for $h(X)$. See also \cite[Proposition 3.2, p.~3]{YY20} where the same arguments are used to verify the conjecture for fibrations which are both $\pi$-trivial and totally non-homologous to zero.

If the fibration is not $\pi$-trivial, which usually is the case, it is in particular necessary to understand how much rational homotopy is contracted when passing to $\pi_*(X)\otimes \qq$. This is dealt with in Theorems \ref{theoB} (in the general situation) and \ref{theoC} (for positively elliptic fibres). So the situation of $\pi$-trivial fibrations, or, more generally, both degeneracy properties of ``$\pi$-triviality'' and ``totally non-homologous to zero", nicely motivate these further generalisations.

\subsection{Convergence}

There are several notions possible for defining ``convergence'' of a family of elliptic spaces. Let us start discussing them.
\begin{defi}
Let $\mathcal{X}$ be a family of elliptic spaces. We say that $\mathcal{X}$ has \emph{accumulation point} $c\in \rr\cup \{\infty\}$ if for any $\varepsilon >0$ there exist infinitely many $X\in \mathcal{X}$ with $|h(X)-c|<\varepsilon$. We say that $\mathcal{X}$ \emph{converges} to $c\in \rr$ if $c$ is its only accumulation point.
\end{defi}

\begin{ex}\label{ex01}
\begin{itemize}
\item Clearly, by definition, no finite family $\mathcal{X}$ of elliptic spaces can have accumulation points nor converge. No family with universally bounded cohomology can converge to zero.
\item Every infinite family $\mathcal{X}$ of elliptic spaces has a (possibly infinite) accumulation point. If $\mathcal{X}$ converges to $c$, then so does any infinite subfamily $\mathcal{Y}\In \mathcal{X}$.
\item The family $\{\cc\pp^n\}_{n\geq 1}$ is a family of universally bounded homotopy $\dim \pi_\odd (\cc\pp^n)\otimes \qq=2$ converging to zero. Compare Proposition \ref{prop03}.
\item The family $\{\s^n\}_{n\geq 2}$ realises infinitely many rational homotopy groups, each element satisfies $\dim \pi_*(X)\leq 2$ for $X\in \mathcal{X}$. Clearly $h(X)\in \{1/2,1\}$ for $X$ in $\mathcal{X}$, and the family has two accumulation points (although the set $\{h(X)\mid X\in \mathcal{X}\}$ is finite and hence does not have any accumulation points). Odd spheres converge to $\tfrac{1}{2}$, even ones to $1$.
\item There are infinite families $\mathcal{X}$ of elliptic spaces realising only finitely many Betti numbers (for example, see \cite[Chapter 6.2, p.~243]{FOT08}). Hence these families have positive accumulation points.
\item Such families can already be found to realise the same cohomology algebras (see \cite{NSY04}).
\item All of these example families only realise finitely many rational homotopy groups, hence the accumulation points are positive, but not infinite. Taking product or more elaborate constructions one may easily adapt limit points.
\end{itemize}
\end{ex}

The next Proposition generalises our observation on the family of spheres.
\begin{prop}
Let $\mathcal{X}$ be a family of elliptic spaces, let $\mathcal{P}$ denote the family of pure spaces, $\mathcal{Q}$ the one of two-stage spaces. If $\mathcal{P}\In  \mathcal{X}$ respectively $\mathcal{Q}\In \mathcal{X}$, then any number $h(X)$ for $X\in \mathcal{P}$ respectively for $X\in \mathcal{Q}$ is an accumulation point of $\mathcal{X}$.
\end{prop}
\begin{prf}
For every pure respectively two-stage space $X$ we construct an infinite sequence of pure respectively two-stage spaces $X_i$ satisfying $h(X)=h(X_i)$ for all $i\geq 1$. This pureness/two-stage property will be obvious from the construction. This can be done as follows.

Recall that pure spaces are two-stage in particular. Let $(\Lambda (V_0\oplus V_1),\dif)$ be the two-stage decomposition of the minimal model of $X$. We choose the minimal model in its isomorphism class such that we display $V_1$ with minimal possible dimension. Hence, the differential is injective on $V_1$ and differentials have a well-defined degree.
Let $v_1,\ldots, v_k$ be a homogeneous basis of $V_0$ and $v'_1,\ldots, v'_{k'}$ be a homogeneous basis of $V_1$. Up to spatial realisation, it suffices to construct a two-stage minimal model $(\Lambda W,\dif)=(\Lambda (W_0\oplus W_1),\dif)$ of $X_i$---for the sake of simplicity we suppress the index $i$ in the models. For this we construct a homogeneous basis $w_0,\ldots, w_k$ of $W_0$ and $w'_1,\ldots, w'_{k'}$ of $W_1$. The $w_j$ and $w_j'$ will be degree shifts of the corresponding $v_j$, $v'_j$. We extend degrees multiplicatively. Hence it remains to define
\begin{align*}
\deg w_j&:={3^i\cdot j}
\\\deg w_j'&:=3^i\cdot\deg (\dif v_j') -1
\end{align*}
and extend degrees multiplicatively as usual. We write $\dif v_j'=p_j(v_i)$ as a polynomial $p_j$ in the $v_i$, and we denote by $p_j(w_i)$ the corresponding polynomial replacing the $v_i$ by the $w_i$. Hence set
\begin{align*}
\dif w_j&:=0
\\\dif w_j'&:=p_j(w_i)
\end{align*}
which is well-defined by construction. Hence all the $X_i$ are well-defined pure respectively two-stage spaces. They are all mutually distinct due to degrees.

Then all $X_i$ have isomorphic minimal models, however, using isomorphisms \emph{not} respecting the grading. Indeed, by construction, the isomorphism to $(\Lambda V,\dif)$ is induced by the correspondence $v_i\sim w_i$, $v_i'\sim w_i'$. (For this note that due to multiplication with $3^i$ the parity of the basis is preserved.) In particular, $h(X_i)=h(X)$ for all $i$. This proves the result.
\end{prf}
As a consequence the family of pure or two-stage spaces or any family containing them like the family of all elliptic spaces does not converge to any limit point.

Note that the elements in the sequences we constructed in the last proof all had the same rational homotopy groups. It seems more interesting to understand what happens if rational homotopy tends to infinity.

\begin{defi}
A family $\mathcal{X}$ of elliptic spaces has \emph{$\pi$-accumulation point} $c\in \rr\cup\{\infty\}$, if
for all $\varepsilon>0$ there exist $n(\varepsilon)\in \nn$ and infinitely many $X\in \mathcal{X}$ with $\dim \pi_*(X)\otimes \qq\geq n(\varepsilon)$ and $|h(X)-c|<\varepsilon$.

The family \emph{$\pi$-converges} to $c\in \rr\cup\{\infty\}$, i.e.~
\begin{align*}
\lim_{\dim \pi_*(X)\otimes \qq\to \infty} h(X)=c
\end{align*}
if $c$ is the only $\pi$-accumulation point.
\end{defi}

In the folllowing let us discuss zero as an accumulation point. We need some preparatory results first.

We provide an easy and coarse estimate on the dimension of the cohomology of a pure space. Note that the important aspect for us is that this estimate is given purely in terms of degrees and dimensions of rational homotopy groups, since the formal dimension $d$ of an elliptic space can be computed just using this degree information.
\begin{lemma}\label{lemma04}
Let $(\Lambda V,\dif)$ be a pure minimal Sullivan algebra of formal dimension $d$. Denote by $a_1,\ldots, a_k$ the degrees of a homogeneous basis of $V^\even$. Then
\begin{align*}
\dim H(\Lambda V,\dif) \leq 2^{\dim V^\odd-\dim V^\even} \cdot \prod_{1\leq i\leq k} \lceil d/a_i \rceil
\end{align*}
\end{lemma}
\begin{prf}
Let $w_1,\ldots, w_k$ be such a homogeneous basis of $V^\even$ with $\deg w_i=a_i$.
Consider the rational fibration given by the relative model
\begin{align}\label{eqn10}
(\Lambda \langle v_1, \ldots, v_k\rangle \otimes \Lambda V,\dif)
\end{align}
with fibre $(\Lambda \langle v_1, \ldots, v_k\rangle,0)$ generated by elements of degrees $\deg v_i=a_i\cdot (\lceil d/a_i\rceil +1)-1$ satisfying $\dif v_i=w_i^{\lceil d/a_i\rceil +1} $.

By construction---we chose the $v_i$ to map to elements of degree larger than the formal dimension $d$ of $(\Lambda V,\dif)$ under the differential $\dif$---the total space actually has the following minimal model up to isomorphism.
\begin{align*}
(\Lambda \langle v_1, \ldots, v_k\rangle \otimes \Lambda V,\dif)
\cong (\Lambda V,\dif)\otimes (\Lambda \langle v_1,\ldots, v_k \rangle,0)
\end{align*}
Hence the relative model \eqref{eqn10} admits a second rational fibration structure with base space $(\Lambda (V^\even\oplus \langle v_1,\ldots, v_k\rangle),\dif)$ and fibre $(\Lambda V^\odd,0)$. The dimension of the cohomology of the base space is $\prod_{1\leq i\leq k} \lceil d/a_i \rceil $.

By the $E_2$-term of the associated Serre spectral sequence of this new fibration we deduce that
\begin{align*}
&\dim H (\Lambda V,\dif) \cdot 2^k
\\=&\dim H((\Lambda V,\dif)\otimes (\Lambda \langle v_1,\ldots, v_k \rangle,0) )
\\  \leq  & \dim H(\Lambda V^\odd,0) \cdot \dim H(\Lambda (V^\even\oplus \langle v_1,\ldots, v_k\rangle),\dif)
\\ \leq & 2^{\dim V^\odd} \cdot \prod_{1\leq i\leq k} \lceil d/a_i \rceil
\end{align*}
The assertion follows.
\end{prf}

For the next proposition it would have been enough to work with the well-known estimate $\dim H^*(X) \leq 2^{\dim X}$ for an elliptic space (and again to use that formal dimension can be expressed via the degrees of rational homotopy groups). As a service to the reader we provided the last lemma with its concise proof instead.
\begin{prop}\label{prop03}
If the family $\mathcal{X}$ of elliptic spaces has $0$ as an accumulation point, then
\begin{align*}
\sup_{X \in \mathcal{X}}\{\dim \pi_*(X)\otimes \qq \}=\infty
\qquad \textrm{or} \qquad
\sup_{ X \in \mathcal{X}}\{i \in \nn \mid \pi_i(X)\otimes \qq \neq 0\}=\infty
\end{align*}

In any case formal dimensions are unbounded, i.e.~
\begin{align*}
\sup_{X \in \mathcal{X}}\{\dim X \}=\infty
\end{align*}
\end{prop}
\begin{prf}
It suffices to show that fixing the rational homotopy groups $\pi_*(X)\otimes \qq$ of an elliptic space $X$, there exists $\alpha\in \nn$ such that $\dim H^*(X')\leq \alpha$ for all elliptic $X'$ satisfying $\pi_*(X')\otimes \qq=\pi_*(X)\otimes \qq$. Indeed, this implies that for $\dim H^*(X)$ to be unbounded within $\mathcal{X}$ (which is clearly necessary for accumulation point zero), it is required to have infinitely many configurations $\pi_*(X)$. That is, either infinitely many homotopy Betti numbers or infinitely many degrees of rational homotopy groups (or both).

In any case the dimension formula (see Section \ref{secrat}) for elliptic spaces (together with the observation that the existence of an even-degree basis element of the rational homotopy groups requires the existence of an additional odd-degree one of at least twice the degree, see \cite[Proposition 32.9]{FHT01})
yields the unboundedness of formal dimensions.

So let us show the existence of $\alpha$. By the odd spectral sequence (\cite[Chapter 32(b), p.~438]{FHT01}) $\dim H(\Lambda V,\dif)\leq \dim H(\Lambda V,\dif_\sigma)$ for a minimal Sullivan algebra $(\Lambda V,\dif)$ with associated pure one $(\Lambda V,\dif_\sigma)$. Hence, without restriction, we may assume that the $X'$ are pure spaces, and we have to show that fixing rational homotopy groups there are only finitely many $\dim H^*(X')$ for pure $X'$ realising the homotopy groups. This follows from Lemma \ref{lemma04} in which we provide an upper bound on cohomology merely in terms of the degrees and dimensions of the rational homotopy groups.
\end{prf}
Recall the family of complex projective spaces from Example \ref{ex01} with constant dimension of rational homotopy groups and diverging cohomology. Here the top degrees of rational homotopy diverge. The family of products of spheres of a fixed dimension clearly has bounded top homotopical degree and diverging homotopical dimension. Both families $\pi$-converge to zero. This illustrates that both cases in the proposition really can occur.

So we already started to answer
\begin{ques*}
Which accumulation points can be realised by a family of elliptic spaces? Or, much more interestingly, which $\pi$-accumulation-points can be realised?
\end{ques*}

\begin{rem}
Instead of merely looking at limits, we also suggest to have a closer look at the rate of convergence.
For example, if convergence is governed by $n\mapsto \tfrac{n}{2^{n}}$, then the elements of $\mathcal{X}$ satisfy the toral rank conjecture. (Of course, this is a rather restrictive condition.)

Clearly, it is well-known (see \cite[Theorem 7.13, p.~279]{FOT08}) that the toral rank $\rk(X)$ of an elliptic $X$ satisfies $\rk(X)\leq \chi_\pi(X)\leq \dim \pi_\odd(X)$. Then
\begin{align*}
\frac{\rk(X)}{\dim H^*(X)}&\leq \frac{\dim \pi_*(X)\otimes \qq}{\dim H^*(X)}\leq \frac{\dim \pi_*(X)\otimes \qq}{2^{\dim \pi_*(X)\otimes \qq}}
\intertext{and (ignoring trivial cases of contractible $X$ or vanishing rank)}
\dim H^*(X) &\geq 2^{\rk(X)} \cdot \frac{ \rk(X)}{\dim \pi_*(X)\otimes \qq}\cdot \frac{2^{\dim \pi_*(X)\otimes \qq}}{2^{\rk(X)}}
\end{align*}
Since $\dim \pi_*(X)\otimes \qq - \rk(X) + \log_2 \rk(X) - \log_2 \dim \pi_*(X)\otimes \qq$ is greater equal $0$ for all relevant values, the toral rank conjecture holds in this situation.
\end{rem}

Next we investigate how convergence to zero behaves under fibrations whence extending the class to more instances.
\begin{prop}\label{prop04}
Let $\mathcal{FIB}=(F\hto{} X \to B)$ be a family of totally non-homologous to zero fibrations of elliptic spaces.
\begin{itemize}
\item
The family $\mathcal{X}$ of total spaces $\pi$-converges to zero if both the family $\mathcal{F}$ of fibres and the family $\mathcal{B}$ of base spaces do.
\item
The family of total spaces $\mathcal{X}$ is $\pi$-convergent to zero if and only if so is the family of spaces $\mathcal{F}\times \mathcal{B}=(F\times B)$.
\end{itemize}
\end{prop}
\begin{prf}
Assume first that both $\mathcal{F}$ and $\mathcal{B}$ $\pi$-converge to zero, and we shall show that $\mathcal{X}$ also $\pi$-converges to $0$. Now with $h(F)$ and $h(B)$ also
\begin{align*}
h(X)=\frac{\dim \pi_*(X)\otimes \qq}{\dim H^*(X)}\leq \frac{\dim \pi_*(F)\otimes \qq+\dim \pi_*(B)\otimes \qq}{\dim H^*(X)}\leq h(F)+h(B)
\end{align*}
tends to zero (using the assumption that both $\dim H^*(F), \dim H^*(B)\leq \dim H^*(X)$).

Due to Lemma \ref{lemma01} and the formula
\begin{align*}
3\dim \pi_*(X)\otimes \qq\geq \dim \pi_*(F)\otimes \qq + \dim \pi_*(B)\otimes \qq
\end{align*}
which we obtain from there, we derive that with $\pi_*(F)\otimes \qq$ and $\pi_*(B)\otimes \qq$ also $\dim \pi_*(X)\otimes \qq$ is unbounded. Hence $\mathcal{X}$ $\pi$-converges to zero.

\bigskip

Now we deal with the second assertion. By the very last argument we derive that $\mathcal{X}$ has unbounded rational homotopy if and only if $\mathcal{F}\times \mathcal{B}$ has. It remains to show that $h(X)$ tends to zero if and only if $h(F\times B)$ does. Due to the fibrations being totally non-homologous to zero we have
\begin{align*}
\tfrac{1}{3}\cdot h(F\times B)&=\frac{\tfrac{1}{3}\cdot \dim \pi_*(F\times B)\otimes \qq}{\dim H^*(X)}
\\\leq h(X)&=\frac{\dim \pi_*(X)\otimes \qq}{\dim H^*(X)}
\\ &\leq  \frac{\dim \pi_*(F)\otimes \qq + \dim \pi_*(B)\otimes \qq}{\dim H^*(F\times B)} =h(F\times B)
\end{align*}
and the assertion follows.
\end{prf}
To avoid confusion: for the family $\mathcal{F}\times \mathcal{B}$ the respective spaces $F$ and $B$ belong to the \emph{same} fibration. As the proof shows, for the first part of the statement instead of a totally non-homologous to zero fibration it would be enough to have the weaker properties $\dim H^*(F), \dim H^*(B)\leq \dim H^*(X)$.


\section{Proof of Theorem \ref{theoA}}\label{secA}

We use the two-stage decomposition for the space $X$ described in the introduction.
\begin{proof}[\textsc{Proof of Theorem \ref{theoA}}]
We need to show that for arbitrarily large numbers there exist infinitely many two-stage spaces $X$ with larger homotopy and with $h(X)$ tending to zero.
It is clear (for example just by taking products of two-stage spaces) that the class of two-stage spaces has unbounded rational homotopy. Hence it remains to see that the number $h(X)$ tends to zero with $\dim \pi_*(X)\otimes \qq$ going to infinity.

We recall from \cite[Theorem 2.3, p.~195]{JL04} that
\begin{align*}
\dim H^*(X)\geq 2^{\dim W^1-\dim V^\even}
\end{align*}
Moreover, by word-length
\begin{align*}
\dim H^*(X)\geq & 1+\dim \Lambda^{\leq 2}V^\even+\dim \Lambda^{\leq 2} W^0\\&+ \dim V^\even\cdot \dim W^0- \dim W^1
\\  =& 1+2\dim V^\even + {\dim V^\even \choose 2} +\dim W^0 + {\dim W^0 \choose 2}\\& + \dim V^\even\cdot \dim W^0 - \dim W^1
\end{align*}
Set $n:=\dim V^\even$, $m:=\dim W^0$, $r:=\dim W^1-\dim V^\even$. (It is clear that $r\geq 0$.)

It follows that
\begin{align*}
h(X)\leq \frac{2n+m+r}{\max (\tfrac{1}{2}(n^2+n+m^2+m+2nm+2-2r),2^r)}
\end{align*}
and we have to show that as one of $n,m,r$ goes to infinity, this expression falls below $1/k$ for any $k\in \nn$.

\bigskip

We consider two different cases.

\case{1} Suppose that $r\leq \frac{n+m}{2}$, which implies that
\begin{align*}
h(X)\leq \frac{5n+3m}{n^2+m^2+2nm+2}
\end{align*}
In this case, whenever $n\to \infty$ or $m \to \infty$ the right hand side becomes arbitrarily small.

\case{2}
Suppose that $r\geq \frac{n+m}{2}$, i.e.~, in particular, $4r\geq 2n+m$. Then $r$ tends to infinity if so do $n$ or $m$. Moreover,
\begin{align*}
h(X)\leq \frac{2n+m+r}{2^r}\leq \frac{5r}{2^r}
\end{align*}
This converges to $0$ whenever any of $n,m,r$ tend to infinity.
\end{proof}
We leave it to the reader to make use of the fact that the estimates are explicit, i.e.~to provide concrete numbers for $n,m,r$ for which the estimates hold.


\section{Proof of Theorem \ref{theoB}}\label{secB}

\begin{rem}\label{remdim}
In the following we shall draw on the dimension formula (see Section \ref{secrat}). We remark that this formula for general elliptic Sullivan algebras $(\Lambda V,\dif)$ is obtained by a reduction to the pure case, i.e.~by passing to the associated pure model $(\Lambda V,\dif_\sigma)$. For this (see \cite[Proposition 32.4, p.~438]{FHT01}) it is shown that in the case when $(\Lambda V,\dif)$ is \emph{minimal} its cohomology is finite dimensional if and only if so is $H(\Lambda V,\dif_\sigma)$. In \cite[Proposition 32.7, p.~442]{FHT01} it is shown that $(\Lambda V,\dif_\sigma)$ and $(\Lambda V,\dif)$ have the same maximal degrees of non-vanishing cohomology. Although, as it seems, despite not being explicitely required in the assertions the proof of this latter result also assumes the minimality of $(\Lambda V,\dif)$. Clearly, already the cohomological finiteness result is wrong without the minimality assumption, as already the example of the contractible algebra $(\Lambda \langle x,y\rangle,x\mapsto y, \deg x=2, \deg y=3)$ with associated pure algebra $(\Lambda \langle x,y\rangle,0)$ of infinite-dimensional cohomology shows. Also \cite[Theorems 32.6, p.~441, and 32.9, p.~442]{FHT01} draw on minimality although, putatively, not stated.

Clearly, the difference between minimal and non-minimal models is eradicated when formulating the dimension formula in terms of homotopy groups (see \cite[p.~434]{FHT01}). The dimension formula, however, stays correct the way it is formulated via even and odd exponents of Sullivan algebras (see Section \ref{secrat}) also for \emph{non-minimal} algebras if either $(\Lambda V,\dif)$ is pure or under the following restriction: Up to isomorphism a Sullivan algebra can be written as the product of a minimal and a contractible one (see \cite[Theorem 14.9, p.~187]{FHT01}). Hence it remains to verify when the dimension formula holds for the contractible factor, i.e.~basically for the two situations $(\Lambda \langle x,y\rangle, x\mapsto y)$ once for $\deg x$ odd and $\deg y=\deg x+1$ even, and once for $\deg x$ even and $\deg y=\deg x+1$ odd. In the first case we obtain dimension $\deg x - (\deg y-1)=\deg x - \deg x=0$, the dimension formula holds; in the second one it fails due to $\deg y - (\deg x-1)=\deg y - \deg y +2=2$. (Note that a pureness assumption excludes the second case.)

However, in the proof of the following lemma we see that the latter algebra cannot be decomposed as the total space of a fibration of elliptic spaces by exactly comparing the dimension formula of such total spaces with the ones of the corresponding product spaces of potential fibre and base. Indeed, this boils down to exactly the same ``$(+2)$-contradiction'' we just observed.
\end{rem}

We now prove a crucial lemma underlying several results. Note that we already drew on it in Section \ref{sec01} (which we do \emph{not} use at all for the subsequent reasoning).
\begin{lemma}\label{lemma01}
Let $F\hto{} X\to B$ be a fibration of rationally elliptic spaces.
Then it holds that
\begin{align*}
\dim \pi_\odd(X)\otimes \qq &\geq \dim \pi_\odd(B)\otimes \qq \geq \dim \pi_\even(B)\otimes \qq \\
\dim \pi_\odd(X)\otimes \qq & \geq \dim \pi_\even(X)\otimes \qq \geq \dim \pi_\even(F)\otimes \qq \\
\dim \pi_\odd(X)\otimes \qq &\geq \dim \pi_\odd (F)
\end{align*}
and, in total,
\begin{align*}
\dim \pi_*(X)\otimes \qq+2\dim\pi_\odd(X)\otimes \qq &\geq \dim \pi_*(F)\otimes \qq + \dim \pi_*(B)\otimes \qq
\end{align*}
\end{lemma}
\begin{prf}
We recall the dimension formula (see Section \ref{secrat}) for the rationally elliptic space $X$.
\begin{align*}
\dim X=\sum_i b_i -\sum_j (a_j-1)
\end{align*}
where the $b_i$ range over the degrees of a homogeneous basis of $\pi_\odd(X)\otimes \qq$, and the analog for the $a_j$ and $\pi_\even(X)\otimes \qq$, and $\dim X$ denotes formal dimension.

We now fix models and homogeneous bases of base and fibre, namely $(\Lambda \langle f_i\rangle_i,\bar \dif)$ a minimal model of $F$, and $(\Lambda \langle b_j\rangle_j,\dif)$ one of $B$ yielding the model of the fibration, i.e.~a (not necessarily minimal) Sullivan model for $X$ given by
\begin{align*}
(\Lambda \langle f_i,b_j\rangle_{i,j},\dif)
\end{align*}
Consider the long exact homotopy sequence
\begin{align*}
\pi_i(F)\otimes \qq\to \pi_i(X)\otimes \qq \to \pi_i(B)\otimes \qq \xto{\partial} \pi_{i-1}(F)\otimes \qq  \to \pi_{i-1}(X)\otimes \qq
\end{align*}
Up to a change of basis we may assume that (passing to the dual sequence) $\ker \partial^*=\langle f_i\rangle_{1\leq i\leq m}$. Consequently (see Section \ref{secrat} and the description of the differential there),
\begin{align*}
\dif_0|_{\langle f_i\rangle_{1\leq i\leq m}}\co \langle f_i \rangle_{1\leq i\leq m} \to \Lambda \langle b_i,f_j\rangle_{i,j} / \Lambda^{\geq 2} \langle b_i,f_j\rangle_{i,j} \cong \langle b_i, f_j\rangle_{i,j}
\end{align*}
is injective with image in $\langle b_i\rangle_i$. Again, up to change of basis, we may assume that $\dif_0(f_i)=b_i$, and $\deg f_i+1=\deg b_i$ for $1\leq i\leq m$. Hence a minimal model of $X$ is given by $(\Lambda \langle f_i, b_j\rangle_{i,j>m},\tilde \dif)$ with a suitably adapted differential $\tilde \dif$.

Next, we use the equality of formal dimensions $\dim X=\dim F + \dim B$ (which can easily be deduced from the Serre spectral sequence and the fact that $E_2^{\dim B,\dim F}=\qq$, which is left invariant by the differentials), and compute both sides separately. By applying the dimension formula to the two respective minimal models of $X$ and of $F\times B$ it follows that
\begin{align*}
&\sum_{i>m} \deg b_i^\odd+ \deg  f_i^\odd -\sum_{j>m} (\deg  b_j^\even+ \deg f_j^\even -2)\\= & \sum_{i} \deg b_i^\odd+\deg f_i^\odd -\sum_{j} (\deg b_j^\even+ \deg f_j^\even -2)
\end{align*}
That is,
\begin{align*}
0&=\sum_{i\leq m} \deg b_i^\odd+\deg f_i^\odd -\sum_{j\leq m} (\deg b_j^\even+ \deg f_j^\even -2)
\\&=\sum_{i\leq m} (\deg f_i^\even+1)- (\deg f_i^\even -1)+\sum_{i\leq m}\deg f_i^\odd -\big( \deg f_i^\odd +1 -1\big)
\\&=2 \cdot \#_{1\leq i\leq m} f_i^\even
\end{align*}
It follows that there is no even-degree element in the kernel of $\dif_0=\partial^*$, i.e.~any even-degree rational homotopy group of $F$ passes non-trivially to $X$. Respectively, the equation writes as
\begin{align*}
0&=\sum_{i\leq m} \deg b_i^\odd+\deg f_i^\odd -\sum_{j\leq m} (\deg b_j^\even+ \deg f_j^\even -2)
\\&=\sum_{i\leq m} \deg b_i^\odd- (\deg b_i^\odd-1 -1)+\sum_{i\leq m}(\deg b_i^\even-1) -\big( \deg b_i^\even -1\big)
\\&=2 \cdot \#_{1\leq i\leq m} b_i^\odd
\end{align*}
and odd-degree rational homotopy groups of the base space $B$ pass injectively to $X$. Both observations taken together prove that
\begin{align*}
\dim \pi_\even(X)\otimes \qq &\geq \dim \pi_\even(F)\otimes \qq \qquad \textrm{and}\\
\dim \pi_\odd(X)\otimes \qq &\geq \dim \pi_\odd(B)\otimes \qq
\end{align*}
It is well known (see \cite[Proposition 32.10, p.~444]{FHT01}) that a rationally elliptic $Y$ satisfies $\dim \pi_\odd(Y)\otimes \qq \geq \dim \pi_\even(Y)\otimes \qq$. Hence it remains to prove that $\dim \pi_\odd(X)\otimes \qq \geq \dim \pi_\odd (F)\otimes \qq$ whence the formula
\begin{align*}
\dim \pi_*(X)\otimes \qq+2\dim\pi_\odd(X)\otimes \qq &\geq \dim \pi_*(F)\otimes \qq + \dim \pi_*(B)\otimes \qq
\end{align*}
follows by summation.

The linear part of the differential $\dif$, namely
\begin{align*}
\dif_0\co \langle f_i^\odd\rangle_i \to \big(\Lambda \langle f_i, b_j\rangle_{i,j}/\Lambda^{\geq 2} \langle f_i, b_j\rangle_{i,j}\big)^\even\cong \langle f_i^\even, b_j^\even\rangle_{i,j}
\end{align*}
maps into $\langle b_i^\even\rangle_i$. From the proof on \cite[p.~443]{FHT01} we cite that for each $b_i^\even$ there exists a basis element $b_i^\odd$ (of degree at least $2\deg b_i^\even -1$). Hence we derive that $\ker \dif_0|_{\langle f_i^\odd\rangle_i}$ passes directly to $\pi_\odd(X)\otimes \qq$, and that also its image, $\im \dif_0$, is injectively represented in the \emph{odd-degree} rational homotopy of $X$. The intersection of those odd degree elements contributed by the fibre and those by the base is clearly trivial. It follows that
\begin{align*}
\dim \pi_\odd (F) \otimes \qq&= \dim \langle f_i^\odd\rangle_i
\\&=\dim \ker \dif_0|_{\langle f_i^\odd\rangle_i}+ \dim \im \dif_0|_{\langle f_i^\odd\rangle_i}
\\&\leq \dim \pi_\odd (X)\otimes \qq
\end{align*}
\end{prf}

\begin{rem}
We remark that the estimate
\begin{align*}
\dim \pi_*(X)\otimes \qq+2\dim\pi_\odd(X)\otimes \qq \geq \dim \pi_*(F)\otimes \qq + \dim \pi_*(B)\otimes \qq
\end{align*}
is sharp as is shown by the example of the Hopf fibration $\s^3\hto{}\s^7\to \s^4$.
\end{rem}

We are finally in the position to provide the
\begin{proof}[\textsc{Proof of Theorem \ref{theoB}}]
From \eqref{eqn123} we recall that $\dim H^*(X)\leq \dim H^*(F)\cdot \dim H^*(B)$. It follows from Lemma \ref{lemma01} that for elliptic spaces $F\hto{} X\to B$ the following estimate holds:
\begin{align}\label{eqn20}
h(F\times B)=\frac{\dim \pi_*(F)\otimes \qq + \dim \pi_*(B)\otimes \qq}{\dim H^*(F)\cdot \dim H^*(B)}\leq  \frac{3\dim \pi_*(X)\otimes \qq}{\dim H^*(X)}=3h(X)
\end{align}
\end{proof}


\section{Proofs of Theorems \ref{theoC} and \ref{theoD}}\label{secCD}

We shall now refine previous arguments to the case of positively elliptic $F$ or $X$.

\begin{proof}[\textsc{Proof of Theorem \ref{theoC}}]
Let us first  prove the right Inequality in \eqref{eqn02} in the depicted cases. For this we observe the following: Due to the multiplicativity of the Euler characteristic in fibrations (see Section \ref{secrat}), the space $X$ is $F_0$ if and only if so are both $F$ and $B$. Hence if $X$ is $F_0$ so are all spaces involved. Moreover, a positively elliptic space has rational cohomology concentrated in even degrees (see \cite[Proposition 32.10]{FHT01}). Hence the Serre spectral sequence degenerates for lacunary reasons, and the fibration is totally non-homologous to zero whence the right inequality in \eqref{eqn02} holds (see Section \ref{secfib}).

The degeneration at the $E_2$-term is enforced by the assumption that $F$ satisfies the Halperin conjecture.

\bigskip

Let us now deal with the left inequality in \eqref{eqn02}. We observed that in both settings from the assertion $F$ is positively elliptic. As in the proof of Theorem \ref{theoB}, we recall that $\dim H^*(X)\leq \dim H^*(F)\cdot \dim H^*(B)$.
From Inequality \eqref{eqn20} we recall that $h(F\times B)\leq 3 h(X)$. In the case when $F$ is positively elliptic we improve this to $h(F\times B)\leq 2 h(X)$ by refining the respective proof. Indeed, it now suffices to show that
\begin{align}\label{eqn70}
2\dim \pi_*(X)\otimes \qq\geq \dim \pi_*(F)\otimes \qq + \dim \pi_*(B)\otimes \qq
\end{align}
since then
\begin{align*}
h(F\times B)&=
\frac{\dim \pi_*(F)\otimes \qq+ \dim \pi_*(B)\otimes \qq}{\dim H^*(F\times B)}
\\&\leq 2\cdot \frac{\tfrac{1}{2}\cdot (\dim \pi_*(F)\otimes \qq + \dim \pi_*(B)\otimes \qq)}{\dim H^*(X)}
\\&\leq 2\cdot \frac{\dim \pi_*(X)\otimes \qq}{\dim H^*(X)}
\\&=2\cdot h(X)
\end{align*}
(Note that the first inequality is actually an equality using that our fibration is totally non-homologous to zero; yet, this is irrelevant for the argument at this stage of the proof.)

As we observed in the proof of Lemma \ref{lemma01}, $(\im \dif_0)_\odd=0$, i.e.~only odd degree homotopy groups from $F$ contract even degree ones from $B$. This implies that
\begin{align*}
\dim \pi_*(X)\otimes \qq=\dim \pi_*(F)\otimes \qq + \dim \pi_*(B)\otimes \qq - 2 c
\end{align*}
where
\begin{align*}
c\leq \min\{\dim \pi_\odd(F)\otimes \qq, \dim \pi_\even(B)\otimes \qq\}
\end{align*}

Since both $F$ is an $F_0$-space, we derive that
\begin{align*}
\dim \pi_\odd(F)\otimes \qq &= \dim \pi_\even(F)\otimes \qq \qquad  \dim \pi_*(F)\otimes \qq=2 \dim \pi_\odd(F)
\end{align*}
Since $\chi_\pi(B)\geq 0$, we always have for elliptic $B$ that
\begin{align*}
\dim \pi_\odd(B)\otimes \qq &\geq \dim \pi_\even(B)\otimes \qq \qquad  \dim \pi_*(B)\otimes \qq\geq 2 \dim \pi_\even(B) \\
\end{align*}
It follows that
\begin{align*}
&\dim \pi_*(X)\otimes \qq
\\\geq &\dim \pi_*(F)\otimes \qq + \dim \pi_*(B)\otimes \qq - \min\{\dim \pi_*(B)\otimes \qq, \dim \pi_*(F)\otimes \qq\}
\\\geq & \max\{\dim \pi_*(B)\otimes \qq, \dim \pi_*(F)\otimes \qq\}
\end{align*}
whence Inequality \eqref{eqn70}.
\end{proof}
\begin{rem}
The inequality \eqref{eqn70} certainly does not hold when $B$ is positively elliptic (instead of $F$). For this just consider the Hopf fibration $\s^3\hto{} \s^7\to \s^4$ (with corresponding inequality $2<3$) again.
\end{rem}

\vspace{5mm}

In the proof of Theorem \ref{theoC} we came to a point where we had to discuss fibrations of $F_0$-spaces. Those are necessarily totally non-homologous to zero. In the proof of Theorem \ref{theoD} we have to deal with fibrations of formal elliptic spaces. As $F_0$-spaces are formal, this generalises the previous discussion. However, such a fibration is no longer necessarily totally non-homologous to zero, as again the example of the Hopf fibration $\s^3\hto{} \s^7 \to \s^4$ shows already. Hence we shall have to discuss the trade-off of homotopy and cohomology degeneration.

Let $F\hto{} X\to B$ be a fibration of formal elliptic spaces. Due to Proposition \ref{prop05} we know that such a formal elliptic space has the structure of the total space of a totally non-homologous to zero fibration of an $F_0$-space over a product of odd-dimensional spheres.
Hence $X$ admits the following Sullivan model
\begin{align*}
(\Lambda V_F\otimes \Lambda T_F \otimes \Lambda V_B\otimes \Lambda T_B,\dif)
\end{align*}
with $T_F$, $T_B$ concentrated in odd degrees, $V_B^\even=V_B^\odd$, $V_F^\even=V_F^\odd$, $(\Lambda V_F\otimes \Lambda T_F,\bar \dif)$ a model of $F$, $(\Lambda V_B\otimes \Lambda T_B,\dif)$ a model of $B$. Next we prove that whenever we contract an element of $T_F$ cohomology halves at least.
\begin{lemma}\label{lemma03}
\begin{align*}
\dim H^*(X)\leq \frac{\dim H^*(F\times B)}{2^{\dim \im (\dif_0|_{T_F})}}
\end{align*}
\end{lemma}
\begin{prf}
We denote by
\begin{align*}
2a_1, \ldots, 2 a_{\dim V_F^\even}, 2 a_{\dim V_F^\even+1}, \ldots, 2 a_{\dim V_F^\even+V_B^\even}
\end{align*}
the degrees of a homogeneous basis of $ V_F^\even \oplus V_B^\even$, and by
\begin{align*}
&2b_1-1, \ldots, 2 b_{\dim V_F^\odd}-1, 2 b_{\dim V_F^\odd+1}-1,\ldots, 2 b_{\dim V_F^\odd+\dim V_B^\odd}-1,
\\& 2 b_{\dim V_F^\odd+\dim V_B^\odd+1}-1, \ldots, 2 b_{\dim V_F^\odd+\dim V_B^\odd+\dim T_F}-1,
\\&2 b_{\dim V_F^\odd+\dim V_B^\odd+\dim T_F+1}-1, \ldots, 2 b_{\dim V_F^\odd+\dim V_B^\odd+\dim T_F+\dim T_B}-1
\end{align*}
the degrees of a homogeneous basis of $V_F^\odd\oplus V_B^\odd \oplus T_F \oplus T_B$.

Since $X$ is formal as well, we can compute its cohomology using the degrees of the rational homotopy groups of $F$ and $B$. Thus it holds that
\begin{align}\label{eqn03}
\dim H^*(X)=  2^{\dim T_F+\dim T_B} \cdot \prod_{1\leq i\leq \dim V_B^\even+\dim V_F^\even} b_{\pi(i)}/ a_i
\end{align}
for some permutation $\pi$ of $\{1, \ldots, \dim V_F^\odd+V_B^\odd+\dim T_F\}$ in particular satisfying $b_{\pi(i)}\leq b_i$. The first factor comes from the product of odd spheres over which $X$ (being formal elliptic) fibres rationally and in a totally non-homologous to zero manner; actually $\dim T_F+\dim T_B=\chi_\pi(X)$. The right hand side computes the cohomological dimensions of possible positively elliptic fibre parts of this totally non-homologous to zero fibration decomposition of $X$. For this we observe that odd-degree homotopy of this part a priori may come from all of $V_F^\odd\oplus V_B^\odd\oplus T_F$. The degree restrictions for the $b_i$ essentially draw on this factor being positively elliptic, i.e.~a non-trivial relation of lower degree cannot be replaced by one of higher degree whence the one of higher degree must be trivial and yields a free factor of odd degree---indeed, the number of relations in the positively elliptic part equals the number of cohomology generators.

As we need to take into account that some homotopy groups might be contracted, we may even have that $b_{\pi(i)}=a_i$. We shall make this more precise. Let $c:=\dim \im (\dif_0|_{T_F})$ denote the dimension of the subgroup of $\pi_\even(B)\otimes \qq$ which is contracted by the rational homotopy groups of $F$ dual to $T_F$ and hence does not contribute to $\pi_\even(X)\otimes \qq$. (We focus on this homotopy solely, although, clearly, more homotopy groups may be contracted by means of $\dif_0(V_F^\odd)$.) We may express the $F_0$-part in the previous estimate as
\begin{align*}
&\prod_{1\leq i\leq \dim V_B^\even+\dim V_F^\even} b_{\pi(i)}/ a_i
\\=&\prod_{1\leq i\leq \dim V_B^\even+\dim V_F^\even-c} \underbrace{b_{\pi(i)}/ a_i}_{\geq 2 \textrm{ or }=1}  \\&\cdot \prod_{\dim V_B^\even+\dim V_F^\even-c+1\leq i\leq \dim V_B^\even+\dim V_F^\even} \underbrace{ b_{\pi(i)}/ a_i}_{=1}
\end{align*}
where we reordered such that the last $c$ factors are those contracted by $T_F$ as depicted. For this, recall again from \cite[p.~443]{FHT01} that, after decomposing the algebra into a minimal one times a contractible one, up to reordering the quotients $b_{\pi(i)}/ a_i$ are at least $2$ on the minimal factor; they equal $1$ on the contractible factor, since, from the proof of Lemma \ref{lemma01} we recall that $\dif_0$ is trivial on $V_F^\even$, and only an odd-degree element of $V_F$ can map non-trivially to $V_B$. Hence the dimension formula yields $b_{\pi(i)}/ a_i=1$ in this situation.

Next we draw some consequences from this description: The minimal model of $F\times B$ is just the product of the minimal models of $F$ and $B$. Hence, in order to compute its cohomology, every factor $b_i/a_i\geq 2$ for $1\leq i\leq \dim V_B^\even+\dim V_F^\even$ yields a factor of at least $2$. In other words, every basis element of $T_F$ contracted via $\dif_0$ hence reduces the dimension of the cohomology of $F\times B$ by a factor of $2$ at least. (Note that this is not true for elements contracted by $\dif_0(V_F^\odd)$.) This together with $b_{\pi(i)}\leq b_i$ from above implies that
\begin{align*}
\dim H^*(X)
 \leq &  2^{\dim T_F+\dim T_B} \cdot \prod_{1\leq i\leq \dim V_F^\even+\dim V_B^\even-c} b_{\pi(i)}/ a_i
 \\ \leq &  2^{\dim T_F+\dim T_B-c} \cdot \prod_{1\leq i\leq \dim V_F^\even+\dim V_B^\even} b_{i}/ a_i
\\ =&  2^{-c} \cdot \dim H^*(F\times B)
\end{align*}
which proves the asserted estimate.
\end{prf}

This now enables us to prove Theorem \ref{theoD} in the form of the next two propositions, one for each estimate in \eqref{eqn02}.
\begin{prop}\label{prop01}
Let $F\hto{} X\to B$ be a fibration of formal elliptic spaces. Then $h(F\times B)\leq 2\cdot h(X)$.
\end{prop}
\begin{prf}
We recall from Lemma \ref{lemma01} that
\begin{align*}
\dim \pi_*(F\times B)\otimes \qq
 \leq & 3\dim \pi_*(X)\otimes \qq
\end{align*}

We combine this with Lemma \ref{lemma03} and the notation from its proof (in particular, $c=\dim \im (\dif_0|_{T_F})$) leading to 
\begin{align*}
h(X)&=\frac{\dim \pi_*(X)\otimes \qq}{\dim H^*(X)}
\\&\geq \frac{2^c\cdot \dim \pi_*(X)\otimes \qq}{\dim H^*(F\times B)}
\\& \geq \frac{2^c\cdot \big(\tfrac{1}{3}\cdot (\dim \pi_*(F\times B)\otimes \qq)\big)}{\dim H^*(F\times B)}
\\&=2^{c}\cdot \tfrac{1}{3} \cdot h(F\times B)
\end{align*}
Hence, in order to establish $h(F\times B)\leq 2 h(X)$, it remains to observe that $2^{c+1} \cdot \tfrac{1}{3}\geq 1$ is equivalent to $c\geq 1$ and to discuss the case $c=0$.

If $c=0$, we argue as follows. We decompose the minimal model of $F$ as $(\Lambda V_F\otimes \Lambda T_F,\bar \dif)$ as in the proof of Lemma \ref{lemma03}. Hence by the arguments from the proof of Lemma \ref{lemma01}, $\dif_0$ can only be non-trivial on a space of dimension $\dim V_F^\odd$. Clearly, $\dim V_F^\odd \leq \tfrac{1}{2} \dim (V_F\oplus T_F)$. Analogously, $\im \dif_0\In V_B^\even$, and $\dim \im \dif_0\leq \tfrac{1}{2}\dim (V_B\oplus T_B)$.
That is, both from fibre and from base space at most half-dimensional rational homotopy is contracted. Hence
\begin{align*}
\dim \pi_*(X)\otimes \qq \geq \frac{1}{2} \dim \pi_*(F\times B)
\end{align*}
Adapting the inequalities above and using \eqref{eqn123}, we then have
\begin{align*}
h(X)&=\frac{\dim \pi_*(X)\otimes \qq}{\dim H^*(X)}
\\& \geq \frac{\tfrac{1}{2}\cdot \dim \pi_*(F\times B)\otimes \qq}{\dim H^*(F\times B)}
\\&= \tfrac{1}{2} \cdot h(F\times B)
\end{align*}
and the result follows also in this case.
\end{prf}

\begin{prop}\label{prop02}
Let $F\hto{} X\to B$ be a fibration of formal elliptic spaces. Then
\begin{align*}
h(X)< h(F)+h(B) + \frac{1}{4}
\end{align*}
\end{prop}
\begin{prf}
The proof basically consists of refining Equation \eqref{eqn03}---in particular, drawing on the terminology and results established there. Hence we recall that
\begin{align*}
\dim H^*(X)=  2^{\dim T_F+\dim T_B} \cdot \prod_{1\leq i\leq \dim V_B^\even+\dim V_F^\even} b_{\pi(i)}/ a_i
\end{align*}
for some permutation $\pi$ of $\{1, \ldots, \dim V_F^\odd+V_B^\odd+\dim T_F\}$ satisfying $b_{\pi(i)}\leq b_i$.

We now claim and prove that up to renumbering $(\pi(1), \ldots, \pi(\dim V_F^\even))=(1,\ldots, \dim V_F^\even)$, i.e.~the first $\dim V_F^\even$ many $b_i$ come from the $F_0$-part of the fibre $F$, i.e.~they are given by the fact that the $2b_i-1$ are the degrees of a homogeneous basis of $V_F^\odd$.

In order to prove this, we again draw on the observations from \cite[p.~443]{FHT01} respectively on \cite[Proposition 32.9, p.~442]{FHT01} and its proof. That is, we have seen that $V_F^\even$ corresponds injectively (respecting degrees) to a subspace of $\pi_\even(X)\otimes \qq$. Hence in the (not necessarily minimal) model of the fibration there must exist a homogeneous subspace $S$ of $V_F^\odd \oplus V_B^\odd \oplus T_F\oplus T_B$ of dimension at least $\dim V_F^\even$ with the property that $\bar \dif S\In \Lambda V_F^\even$ and that $(\Lambda V_F^\even\otimes \Lambda S,\dif)$ is elliptic. Here, as usual, $\bar \dif$ denotes the projection of $\dif$ to the fibre $\Lambda (V_F\oplus T_F)$. Since $\bar \dif|_{T_F\oplus V_B\oplus T_B}=0$, the only such subspace is gradedly isomorphic to $V_F^\odd$ itself. In other words, since $\dim V_F^\odd=\dim V_F^\even$, these first degrees $b_i$ for $1\leq i\leq \dim V_F^\even$ are uniquely determined by a homogeneous basis of $V_F^\odd$.

Hence using \eqref{eqn098} we can refine Formula \eqref{eqn03} by
\begin{align}\label{eqn05}
\dim H^*(X)= 2^{\dim T_B} \cdot \dim H^*(F)\cdot \prod_{1\leq i\leq \dim V_B^\even} b_{\pi(i)}/ a_i
\end{align}
where the $b_i$ are the odd exponents of $V_B\oplus T_F$, and $\pi$ now is a permutation of $\{1, \ldots, V_B^\odd+\dim T_F\}$ satisfying $b_{\pi(i)}\leq b_i$. Indeed, we have seen that both even and also odd exponents of the $F_0$-part of $F$ appear in the product. That is, the product of their quotients computes the dimension of the cohomology of the $F_0$-part of $F$, i.e.~$\dim H(\Lambda V_F,\bar \dif)$. With the decomposition of formal elliptic spaces (see Proposition \ref{prop05} and Remark \ref{remmin}) yielding $\dim H^*(F)=\dim H(T_F,0) \cdot \dim H(\Lambda V_F,\bar \dif)$, and with $2^{\dim T_F}=\dim H(T_F)$ the refined formula follows.

We derive the estimate
\begin{align}\label{eqn04}
\nonumber\dim H^*(X)&\geq  2^{\dim T_B+\dim V_B^\even} \cdot \dim H^*(F) \intertext{or, equivalently,}
\dim H^*(X)&\geq  2^{\chi_\pi(B)+\dim \pi_\even(B)\otimes \qq} \cdot \dim H^*(F)
\end{align}

Clearly, we have that $\chi_\pi(B)+\dim \pi_\even(B)\otimes \qq\geq \tfrac{1}{2}\cdot \dim \pi(B)\otimes \qq$, and $\dim \pi_*(X)\otimes \qq\leq \dim \pi_*(F)\otimes \qq +\dim \pi_*(B)\otimes \qq$.

Hence we can estimate
\begin{align}\label{eqn08}
h(X)=&\frac{\dim \pi_*(X)\otimes \qq}{\dim H^*(X)}
\nonumber \\ \leq  &\frac{\dim \pi_*(X)\otimes \qq}{2^{\chi_\pi(B)+\dim \pi_\even(B)\otimes \qq} \cdot \dim H^*(F)}
\nonumber \\ \leq & \frac{\dim \pi_*(F)\otimes \qq+\dim \pi_*(B)\otimes \qq}{2^{(\dim \pi_*(B)\otimes \qq)/2} \cdot \dim H^*(F)}
\nonumber\\= & \frac{\dim \pi_*(F)\otimes \qq}{2^{(\dim \pi_*(B)\otimes \qq)/2} \cdot \dim H^*(F)} + \frac{\dim \pi_*(B)\otimes \qq}{2^{(\dim \pi_*(B)\otimes \qq)/2} \cdot \dim H^*(F)}
\end{align}
The formula $h(X)< h(F)+ h(B) + \tfrac{1}{4}$ trivially holds true whenever one of $F$ and $B$ are contractible. Hence we may assume this not to be the case. As an elliptic space satisfies Poincar\'e duality, it follows that both $\dim H^*(F), \dim H^*(B) \geq 2$. We derive that
\begin{align}\label{eqn06}
h(X)< h(F) + \frac{\dim \pi_*(B)\otimes \qq}{2^{(\dim \pi_*(B)\otimes \qq)/2+1} }
\end{align}
and we need to discuss the inequality $\tfrac{n}{2^{n/2+1}}\leq \frac{1}{4}$ for $n\in \nn$ (playing the role of $\dim \pi_*(B)\otimes \qq$). This holds true unless $n\in [1,7]$.

\bigskip

So, in order to finish the proof, we need to differ and discuss the following particular cases:
\begin{itemize}
\item[(i)] $\dim H^*(F)=2$, implying that either
\begin{enumerate}
\item $F\simeq \s^{2k+1}$, $k\geq 0$, or
\item $F\simeq \s^{2k}$, $k\geq 1$.
\end{enumerate}
\item[(ii)] $\dim H^*(F)=3$ equivalent to $F\simeq_\qq \qq[x]/x^3$.
\item[(iii)] $\dim H^*(F)\geq 4$
\end{itemize}

\bigskip

\case{(i.1)}
Let us first deal with Case (i.1). That is, we consider a fibration
\begin{align*}
\s^1 \hto{} X\to B
\end{align*}
with $\dim \pi_*(B)\otimes \qq\leq 7$. It follows that $\dim \pi_*(X)\otimes \qq\leq 8$. As $X$ is formal, from Proposition \ref{prop05} we derive that, depending on the dimension of its rational homotopy, the cohomology of $X$ satisfies the following: $(\dim \pi_*(X)\otimes \qq,H^*(X))$ can be estimated from below by $(n, \geq 2^{\lceil n/2\rceil})$. Correspondingly, in the respective cases, $(\dim \pi_*(X)\otimes \qq ,h(X))$ can be estimated by $(1,\leq \tfrac{1}{2} )$, $(2,\leq 1)$, $(3,\leq \tfrac{3}{4})$, $(4,\leq 1)$, $(5,\leq \tfrac{5}{8})$, $(6,\leq \tfrac{3}{4})$, $(7,\leq \tfrac{7}{16})$, $(8,\leq \tfrac{1}{2})$.

Since in our case $h(F)=h(\s^{2k+1})=1/2$, we derive that the inequality $h(X)\leq 1/2 + 1/4=3/4$---and hence the asserted strict inequality $h(X)<h(F)+h(B)+1/4$---holds unless $\dim \pi_*(X)\otimes \qq\in  \{2,4\}$. In these latter two cases the estimates, however, are sharp only when $X$ is positively elliptic. That is, given that by the additivity of the homotopy Euler characteristic $\chi_\pi(X)\geq \chi_\pi(\s^{2k+1})=1$, the actual upper bounds in these two cases are given by $(2,\geq 4)$, $(4,\geq 16)$ for $(\dim \pi_*(X)\otimes \qq,H^*(X))$ and by $(2,\leq \tfrac{1}{2})$, $(4,\leq \tfrac{1}{4})$ for $(\dim \pi_*(X)\otimes \qq,h(X))$. Hence we are also done in these cases.

We remark that the arguments underlying this are our usual estimates of the cohomology of the $F_0$-factor: given the decomposition $(\Lambda B\otimes \Lambda V,\dif)$ from Proposition \ref{prop05} and Remark \ref{remmin} we estimate $\dim H(\Lambda B,0)\geq 2^{\dim B}$ and $\dim H(\Lambda V,\dif)\geq 2^{\dim V/2}$. That is, once we have fixed the dimension $\dim \pi_*(X)=\dim V+\dim B$ (see Remark \ref{remmin}) of the total rational homotopy, the smaller the dimension of the rational homotopy of the $F_0$-part, $\dim V$, the larger the overall cohomology $\dim H^*(X)$ predicted by this estimate.

\bigskip

\case{(i.2)} \case{(ii)}
Since the Halperin conjecture is confirmed for spaces with cohomology algebra generated by one element (see Section \ref{secrat}; for monicly generated cohomology this is just a trivial computation), the fibration is totally non-homologous to zero in Cases (i.2) and (ii). As we recalled in Section \ref{secfib}, the formula $h(X)\leq h(F) + h(B)$ holds whenever the fibration is totally non-homologous to zero. Hence we are done in these cases.

\bigskip

\case{(iii)}
In Case (iii) we refine Inequality \eqref{eqn06} to
\begin{align}\label{eqn07}
h(X)< h(F) + \frac{\dim \pi_*(B)\otimes \qq}{2^{\dim \pi_*(B)\otimes \qq/2+2} }
\end{align}
and solve that $\tfrac{n}{2^{n/2+2}}\leq \frac{1}{4}$ holds for $n\in \nn\setminus \{3\}$. Hence, eventually, assume $\dim \pi_*(B)\otimes \qq=n=3$. Let us see that this case is merely an artefact of the proof. Indeed, we provide a refined version of Estimate \eqref{eqn08} directly derived using Inequality \eqref{eqn04} and not the simplification $\chi_\pi(B)+\dim \pi_\even(B)\otimes \qq\geq \tfrac{1}{2}\cdot \dim \pi(B)\otimes \qq$. That is, we obtain
\begin{align*}
h(X)< h(F) + \frac{\dim \pi_*(B)\otimes \qq}{2^{\chi_\pi(B)+\dim \pi_\even(B)\otimes \qq} \cdot \dim H^*(F)}
\end{align*}
The cases when $\dim \pi_*(B)\otimes \qq=n=3$ now are the following:
\begin{itemize}
\item Either $\chi_\pi(B)=3$ and, due to formality, $B$ rationally is a product of three odd-dimensional spheres, or
\item $\chi_\pi(B)=1$ and $B$ has an $F_0$-component with cohomology algebra generated by one element.
\end{itemize}
In the first case respectively the second case we derive that
\begin{align*}
\frac{\dim \pi_*(B)\otimes \qq}{2^{\chi_\pi(B)+\dim \pi_\even(B)\otimes \qq} \cdot \dim H^*(F)} & \leq \frac{3}{2^{3}\cdot 4}<\frac{1}{4}
\\ \frac{\dim \pi_*(B)\otimes \qq}{2^{\chi_\pi(B)+\dim \pi_\even(B)\otimes \qq} \cdot \dim H^*(F)} & \leq \frac{3}{2^{2}\cdot 4}<\frac{1}{4}
\end{align*}
and we are done.
\end{prf}

As a corollary of the proof let us fix Observation \eqref{eqn04} again, as it may be of independent interest.
\begin{cor}
For a fibration of formal elliptic spaces $F\hto{} X\to B$ we have the estimate
\begin{align*}
\dim H^*(X)&\geq  2^{\chi_\pi(B)+\dim \pi_\even(B)\otimes \qq} \cdot \dim H^*(F)
\end{align*}
\end{cor}
\vproof

We sum up the results of these propositions.
\begin{proof}[\textsc{Proof of Theorem \ref{theoD}}]
Clearly, Theorem \ref{theoD} is a combination of Propositions \ref{prop01}, and \ref{prop02}.
\end{proof}

\vspace{5mm}

We finally prove the Conjecture for an interesting class of manifolds, namely for any fibration with $X$ rationally one of the known simply-connected manifolds of positive sectional curvature. In particular, this includes all simply-connected homogeneous spaces admitting homogeneous metrics of positive curvature---see Section \ref{secpos} for details. The key observation for this is that any such manifold $M$ is formal and has the rational structure of an $F_0$-space, if $\dim M$ is even. If $\dim M$ is odd, they satisfy
\begin{align}\label{typ01}
\tag{$*$}
\dim \pi_\odd(M)\otimes \qq= 2 \qquad \textrm{and} \qquad  \dim \pi_\even(M)\otimes \qq=1
\end{align}
unless $M$ rationally is an odd-dimensional sphere.
\begin{rem}\label{rem01}
Indeed, clearly, the class of spaces with finite dimensional rational cohomology and satisfying $(*)$ is exactly the class of spaces with rational cohomology algebra generated by one even-degree and one odd-degree element.
\end{rem}
\begin{lemma}\label{lemma02}
An elliptic Sullivan algebra $(\Lambda V,\dif)$ satisfying $(*)$ is formal.
\end{lemma}
\begin{prf}
A minimal model of this algebra is of the form \linebreak[4]$(\Lambda V,\dif)=(\Lambda \langle x,y,z\rangle, \dif)$ with $\deg x$ even, $\deg y$, $\deg z$ odd. Since the associated pure algebra has finite-dimensional cohomology if and only if the original one has, we derive that, without restriction, $\deg z\geq \deg y$ and $\deg z>\deg x$. We derive that $(\Lambda V,\dif)$ decomposes as the total space of a fibration with fibre $(\Lambda \langle x,z\rangle,\bar \dif)$ over $(\Lambda \langle y\rangle,0)$. Since $(\Lambda V,\dif)$ is simply-connected whence $\deg y>1$, we obtain that $\dif x=0$. Consequently, $(\Lambda V,\dif)\cong (\Lambda \langle x,z\rangle, \dif)\otimes (\Lambda \langle y\rangle,0)$  with $\dif x=0$ and $\dif z=x^k$ for some $k>0$. Such an algebra is clearly formal.
\end{prf}

\begin{proof}[\textsc{Proof of Corollary \ref{corB}}]
We need to discuss the potential fibrations for all total spaces $X$ which we depicted around $(*)$.

\case{1} If $X$ is even-dimensional, then $X$ is an $F_0$-space. Hence the result follows from Theorem \ref{theoC}.

\case{2} Let us assume that $X$ rationally is an odd-dimensional sphere.
By Lemma \ref{lemma01} we obtain that $\dim \pi_\even(F)\otimes \qq=0$, $\dim \pi_\odd (B)\otimes \qq\leq 1$, whence $\dim \pi_\even(B)\otimes \qq\leq 1$ and $\dim \pi_\odd(F)\otimes \qq\leq 1$. In total, it follows that both $F$ and $B$ have one of the models $(\Lambda \langle x\rangle,0)$, $\deg x$ odd, or $(\Lambda \langle x,y\rangle, \dif)$ with $\deg x$ even, $\deg y$ odd, $\dif x=0$, $\dif y=x^k$ for some $k>1$. Both are formal, and the result follows from Propositions \ref{prop01} and \ref{prop02}.

\case{3} Now suppose $X$ is odd-dimensional and not a sphere. Then $(*)$ applies. By the additivity of the homotopy Euler characteristic we know that $1=\chi_\pi(M)=\chi_\pi(F)+\chi_\pi(B)$. Hence either $\chi_\pi(F)=0$ and $\chi_\pi(B)=1$ or $\chi_\pi(F)=1$ and $\chi_\pi(B)=0$. In the first case $F$ is positively elliptic. Moreover, since $\pi_\even(F)\otimes \qq \leq \dim \pi_\even(M)\otimes \qq$, it follows that $H^*(F)$ is generated by one element. Since the Halperin conjecture is confirmed for at most $3$ cohomology algebra generators (see Section \ref{secrat}), the result follows again from Theorem \ref{theoC}.

In the second case, $\chi_\pi(F)=1$, the space $B$ is positively elliptic, and we shall have to distinguish yet two more non-trivial cases (using Lemma \ref{lemma01} again). For this we first note that $\dim \pi_\even(F)\otimes \qq\leq 1$, $\dim \pi_\odd(F)\otimes \qq\leq 2$, $\dim \pi_\odd(B)\otimes \qq\leq 2$, $\dim \pi_\odd(B)\otimes \qq\leq 2$. Combining these pieces of information leads to the following cases.
\begin{itemize}
\item[(i)]
$F$ either rationally is an odd-dimensional sphere and $\dim \pi_\even(B)\otimes \qq=\dim \pi_\odd(B)=1$, or $\dim \pi_\even(B)\otimes \qq=\dim \pi_\odd(B)=2$, or
\item[(ii)]
$F$ is of the type $(\ast)$ described afore the proof. Hence $H^*(B)$ is either generated by one element, by two elements or contractible.
\end{itemize}
In order to apply Propositions \ref{prop01} and \ref{prop02} it remains to observe using Lemma \ref{lemma02} that in any case all of $X$, $F$, and $B$ are formal.
\end{proof}
As we noted in \cite[Corollary 4.15, p.~2292]{AK15} there are several geometric conditions which require a positively curved manifold to be a compact rank one symmetric space, and hence to satisfy Conjecture \ref{conj01} in particular, like a $2$-positive curvature operator, weakly quarter-pinched curvature, or the sectional curvature bound $\sec\geq 1$ together with the diameter bound $\operatorname{diam} M\geq \pi/2$.



\def\cprime{$'$}


\vfill

\begin{center}
\noindent
\begin{minipage}{\linewidth}
\small \noindent \textsc
{Manuel Amann} \\
\textsc{Institut f\"ur Mathematik}\\
\textsc{Differentialgeometrie}\\
\textsc{Universit\"at Augsburg}\\
\textsc{Universit\"atsstra\ss{}e 14 }\\
\textsc{86159 Augsburg}\\
\textsc{Germany}\\
[1ex]
\textsf{manuel.amann@math.uni-augsburg.de}\\
\end{minipage}
\end{center}

\end{document}